# ON THE PERFORMANCE OF FDR CONTROL: CONSTRAINTS AND A PARTIAL SOLUTION[1]

By Zhiyi Chi

*University of Connecticut*

The False Discovery Rate (FDR) paradigm aims to attain certain control on Type I errors with relatively high power for multiple hypothesis testing. The Benjamini–Hochberg (BH) procedure is a well-known FDR controlling procedure. Under a random effects model, we show that, in general, unlike the FDR, the positive FDR (pFDR) of the BH procedure cannot be controlled at an arbitrarily low level due to the limited evidence provided by the observations to separate false and true nulls. This results in a criticality phenomenon, which is characterized by a transition of the procedure's power from being positive to asymptotically 0 without any reduction in the pFDR, once the target FDR control level is below a positive critical value. To address the constraints on the power and pFDR control imposed by the criticality phenomenon, we propose a procedure which applies BH-type procedures at multiple locations in the domain of $p$-values. Both analysis and simulations show that the proposed procedure can attain substantially improved power and pFDR control.

**1. Introduction.** Since the original work of Benjamini and Hochberg [2], the False Discovery Rate (FDR) paradigm has become an attractive approach to multiple hypothesis testing (cf. [1, 2, 5, 6, 7, 8, 9, 10, 12, 15, 16, 17] and references therein). FDR is the expected value of the false discovery proportion (FDP), that is, the proportion of falsely rejected hypotheses among all those rejected if there is at least one rejection, and 0 otherwise. Denoting by $R$ the total number of rejected nulls, and $V$ that of rejected true nulls, $\text{FDP} = \frac{V}{R \vee 1}$ and $\text{FDR} = E[\text{FDP}]$. One of the most well-known FDR controlling procedures is the Simes procedure [14] adopted by Benjamini and Hochberg [2], henceforth referred to as the BH procedure. It is now a

Received December 2005; revised July 2006.

[1]Supported in part by NIH Grant MH68028.

classical result that, under certain conditions on the $p$-values, for any target control level $\alpha \in (0,1)$, the BH procedure can attain FDR $\leq \alpha$ [2, 7, 17].

One important issue related to FDR control is power. Let $n$ be the number of nulls being tested and $N$ the number of true nulls among them. Then power $= \frac{R-V}{(n-N)\vee 1}$. A main merit that FDR control is thought to have is its relatively high power compared to Familywise Error Rate (FWER) control, which is on $P(V > 0)$. Many papers show that incorporating an estimate of $\frac{N}{n}$ can increase the power of the BH procedure [3, 9, 15, 18].

Another issue related to FDR control is positive FDR (pFDR), which is the expectation of FDP, conditional on there being at least one rejection: $E[\frac{V}{R}|R \geq 1]$ [15, 16]. Conceptually, pFDR is important for follow-up studies once discoveries are made. However, it is known that for FDR control with a fixed rejection region, pFDR $\to$ FDR as $n \to \infty$ [8, 15]. Therefore, pFDR has been practically treated the same as FDR for $n \gg 1$.

Despite the importance of power and pFDR, there seems to be little work on whether there are any constraints on them. We shall take an asymptotic approach to this issue, as FDR control is often applied when $n$, the number of tested hypothesis, is large. There are two basic and interrelated questions. First, is pFDR always asymptotically the same as FDR? Second, can the BH procedure always attain an asymptotically positive power for a target FDR control level?

The relevance of these two questions will be illustrated by several examples under the setting of a random effects model. The examples include $t$-tests, $F$-tests and multiple tests where the false null distribution is a mixture of Gaussians with variances smaller than the true null Gaussian distribution. As will be seen, in each example there is a critical value $\alpha_* > 0$ for the target FDR control level $\alpha$ giving rise to different regimes of behavior of the BH procedure. When $0 < \alpha < \alpha_*$, the power of the BH procedure decays to 0 at the rate $O_p(1/n)$. Meanwhile, the pFDR converges to a certain constant $\beta_*$ which is strictly greater than the FDR. As a result, the pFDR and FDR are different. When $\alpha > \alpha_*$, the power of the BH procedure converges to a positive constant. However, while the pFDR and FDR are asymptotically equal, both are strictly greater than $\beta_*$. Thus, the pFDR is always bounded below by $\beta_*$, which actually can be quite large if the proportion of false nulls is low. In contrast, the FDR can be controlled at any level. This controllability, however, has a cost: when $\alpha < \alpha_*$, the power of the BH procedure is of the same order as the power of a FWER controlling procedure.

Importantly, the above "criticality" phenomenon is not peculiar to the BH procedure. Once test statistics are selected, for all the multiple testing procedures based on them, there is a common, possibly positive lower bound on the pFDR. Pushing the FDR below this bound separates the FDR and pFDR, and leads to asymptotically zero power. The bound is "intrinsic" in that it is purely a consequence of the distributional properties of the test statistics. How the bound affects power and pFDR control for multiple testing in general is studied elsewhere [4].

In view of the criticality phenomenon, it is natural to explore ways to improve the performance of FDR control. As demonstrated by much work, this is possible by appropriately increasing the target FDR control level of the BH procedure [3, 9, 15, 18]. In this article we propose a procedure which applies BH-type procedures at different locations, or "reference points" in the domain of $p$-values. The idea is to utilize the distributional properties of the $p$-values more effectively, which contain important information that distinguishes false nulls from true nulls. When the distribution function of the $p$-values is concave, as in the case for $t$- and $F$-tests, this procedure is asymptotically the same as the BH procedure. However, in general, as in the case involving Gaussian mixtures, the power and pFDR control can be improved significantly. Unlike the BH procedure, which generates a single random rejection interval containing 0, the multi-reference point procedure can generate multiple random rejection intervals.

The rest of this article is organized as follows. Section 2 collects the main theoretical results on the criticality phenomenon of the BH procedure. We identify the critical value for the target FDR control level and the lower bound for pFDR, and state various asymptotics of the BH procedure. Section 3 considers examples involving multiple $t$-, $F$- and $z$-tests. For the first two, the strict concavity of the distribution of $p$-values will be shown. The limitation of the BH procedure for the $z$-tests provides some motivation for Section 4, which proposes the aforementioned multi-reference point procedure. Section 5 reports a numerical study. Section 6 concludes with some remarks. Details of the proofs are given in the supplemental material.

The rest of this section covers preliminaries. Given unadjusted marginal $p$-values $p_1, \ldots, p_n$, each one for a different null hypothesis, let $p_{n:1} \leq \cdots \leq p_{n:n}$ be their order statistics. Set $p_{n:0} = 0$ and $p_{n:n+1} = 1$. Given the target FDR control level $\alpha \in (0,1)$, the BH procedure rejects all hypotheses with $p$-values $\leq p_{n:R_n}$, where

$$(1.1) \qquad R_n = \max\left\{j \geq 0 : p_{n:j} \leq \frac{\alpha j}{n}\right\}.$$

The number of false rejections and power are

$$V_n = \#\{j \leq n : p_j \leq p_{n:R_n},\ \text{the } j\text{th null is true}\}, \qquad \text{power}_n = \frac{R_n - V_n}{(n - N) \vee 1}.$$

In the following sections, the $p$-values are assumed to be sampled from a random effects model as follows. Let the population fraction of false nulls among all the nulls be a fixed $\pi \in (0,1)$. Let $\theta_j := \mathbf{1}\{\text{the } j\text{th null is false}\}$. Then $(p_1, \theta_1), (p_2, \theta_2), \ldots$ are i.i.d., such that $\Pr\{\theta_j = 1\} = \pi$ and

$$\Pr\{p_j \leq u | \theta_j = 0\} = u, \qquad \Pr\{p_j \leq u | \theta_j = 1\} = G(u), \qquad u \in [0,1].$$

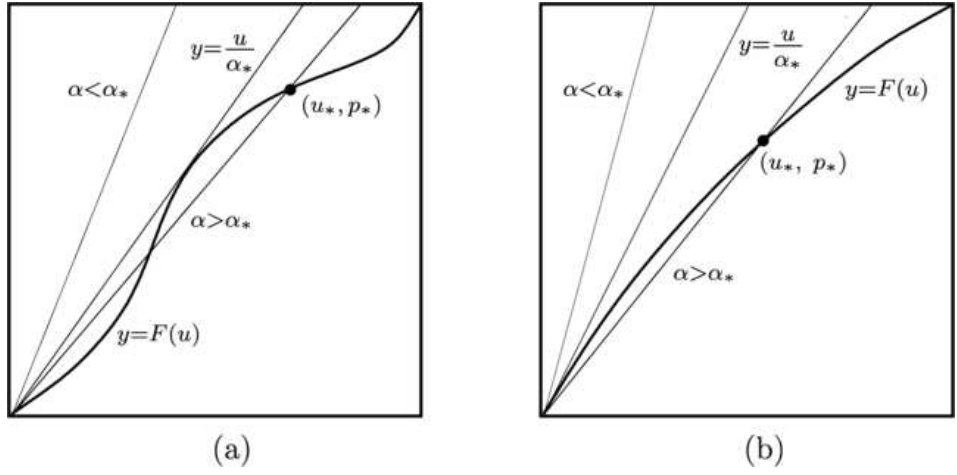

FIG. 1.  *Critical value $\alpha_*$ for general $F$ (a) and concave $F$ (b).*

Under this model, FDR $= (1-\pi)\alpha$ for the BH procedure [2, 7, 17].

**2. Main theoretical results.** This section collects analytical results on the criticality phenomenon of the BH procedure, that is, there can be a critical value $\alpha_* > 0$ such that, for $\alpha < \alpha_*$ and $\alpha > \alpha_*$ the asymptotic behavior of the procedure is categorically different.

By the random effects model, the common distribution function of the $p$-values is

(2.1) $$F(u) = (1-\pi)u + \pi G(u).$$

We shall always assume that $F \in C([0,1])$ with $F(0) = 0$. Denote

(2.2) $$\alpha_* = \inf_{u>0} \frac{u}{F(u)} \leq 1, \qquad \beta_* = (1-\pi)\alpha_*.$$

Figure 1 illustrates the meaning of $\alpha_*$. In particular,

(2.3) $\quad F$ is strictly concave $\implies \alpha_* = \dfrac{1}{F'(0)}, \qquad \beta_* = \dfrac{1-\pi}{F'(0)}.$

To get an idea why $\alpha_*$ is the critical value, observe $R_n = \max\{j : p_{n:j} \leq \alpha F_n(p_{n:j})\}$, where $F_n$ is the empirical distribution of the $p$-values. For $n \gg 1$, $F_n \approx F$. Then the largest rejected $p$-value $\approx u_*$, and $R_n/n \approx p_*$, where

$$u_* = \max\{u \in [0,1] : u/\alpha \leq F(u)\}, \qquad p_* = F(u_*) = u_*/\alpha.$$

To be more precise, in a certain probabilistic sense,

(2.4)
$$u_* = \lim_n [\text{Largest rejected } p\text{-value}],$$
$$p_* = \lim_n [\text{Proportion of rejected } p\text{-values}].$$

From Figure 1, it can be seen that for $\alpha < \alpha_*$, $u_* = p_* = 0$, which suggests that $R_n$ is of order $o(n)$ and its asymptotic behavior should depend only on the local properties of $y = u/\alpha$ and $y = F(u)$ around 0, especially $\alpha$ and $F'(0)$. On the other hand, for $\alpha > \alpha_*$, $p_* > 0$ and hence $R_n \approx np_*$. Finally, if $y = u/\alpha$ and $y = F(u)$ are only tangent at 0 as in Figure 1(b), then a third type of asymptotic behavior arises when $\alpha = \alpha_*$.

The asymptotic distribution of the BH procedure when $\alpha < \alpha_*$ can be characterized as follows. Note that $1 - \pi \leq F'(0) \leq 1/\alpha_*$.

THEOREM 2.1. *Suppose $\alpha \in (0, \alpha_*)$. Let $c = \alpha F'(0)$. Then, as $n \to \infty$,*

$$R_n \xrightarrow{d} \tau := \max\{j : S_j < 0\},$$

*the last time of excursion into $(-\infty, 0)$ by the random walk $S_0 = 0$, $S_j = S_{j-1} + \gamma_j - c$, $j \geq 1$, with $\gamma_1, \gamma_2, \ldots$ i.i.d. $\sim \mathrm{Exp}(1)$ with density $e^{-x}$, $x > 0$. The distribution of $\tau$ is*

$$(2.5) \qquad \Pr\{\tau = k\} = \frac{k^k}{k!}(1-c)c^k e^{-kc}, \qquad k = 0, 1, \ldots.$$

*Consequently, the power of the BH procedure is of order $O_p(1/n)$. Furthermore,*

$$(2.6) \qquad \mathrm{pFDR} \to \beta_*$$

*and*

$$(2.7) \qquad \sum_{k=1}^{\infty} d_{\mathrm{TV}}(\mathcal{L}(V_n | R_n = k), \mathrm{Binomial}(k, \beta_*)) \Pr\{R_n = k\} \to 0,$$

*where $d_{\mathrm{TV}}(\mu, \nu) := \sum_k |\mu_k - \nu_k|$ denotes the total variation distance of two distributions $\mu$ and $\nu$ on $0, 1, 2, \ldots$, and $\mathcal{L}(V_n | R_n = k)$ is the conditional distribution of $V_n$.*

The limiting distribution of $R_n$ was established for $\pi = 0$ in [7]. The characterization of the distribution in terms of an excursion time is new in the context of FDR. The implication of Theorem 2.1 is much stronger than that in [7]. It shows that when $\alpha < \alpha_*$, the total number of true discoveries is bounded, no matter how large $n$ is, even though there is a fixed positive fraction of false nulls. Equation (2.6) indicates that the pFDR cannot be lower than $\beta_* > 0$. Equation (2.7) means that, given a nonzero value of $R_n$, the conditional distribution of $V_n$ approximately follows a Binomial distribution.

Next we consider the case where $\alpha_* \leq \alpha < 1$. To establish the asymptotics, we assume that all the multiple tests are based on an infinite sequence of hypotheses with corresponding $p$-values $p_1, p_2, \ldots$, such that $R_n$ and $V_n$ are attained from the BH procedure applied to the first $n$ of them. We first consider the practically important case $\alpha > \alpha_*$. Then $0 < u_* < p_* < 1$. In [8] it is shown that the proportion of rejected hypotheses $R_n/n \to p_*$ in probability. The law of iterated logarithm (LIL) below characterizes the convergence of $R_n/n$ in a stronger sense, namely, almost sure convergence.

THEOREM 2.2.  *Suppose $\alpha \in (\alpha_*, 1)$ and $\Delta = 1 - \alpha F'(u_*) > 0$. Let $q_* = 1 - p_*$. Then*

$$\limsup_{n} \pm \frac{R_n - np_*}{\sqrt{n \log \log n}} = \frac{\sqrt{2 p_* q_*}}{\Delta}, \qquad a.s. \tag{2.8}$$

*Furthermore, $R_n/n$ is asymptotically proportional to the power:*

$$\text{power}_n = \frac{R_n}{n}\left(\frac{1-\alpha}{\pi} + \alpha\right) + o_p(1) \to G(u_*), \qquad a.s. \tag{2.9}$$

The condition on $\Delta$ means that the line $y = u/\alpha$ crosses the graph of $y = F(u)$ instead of being tangent at $(u_*, p_*)$; see Figure 1. Theorem 2.2 implies that $R_n/n \to p_*$ at rate $O(\sqrt{\log \log n / n})$. Since now the pFDR is asymptotically equal to the FDR, and the latter is equal to $(1-\pi)\alpha$ [2, 7, 17], the pFDR is asymptotically strictly greater then $\beta_*$. If $\pi < 1$, then (2.9) implies

$$\text{power}_n > R_n/n, \qquad \text{for } n \gg 1.$$

The inequality provides a conservative estimate of the power, which can be useful since neither the number of false nulls nor that of rejected false nulls is directly observable.

The main tool to prove Theorem 2.2 is Kiefer's result on Bahadur representation ([13], Section 15.1). Details of the proof are given in the supplemental materials.

Conceptually, it is of interest to consider the behavior of the BH procedure when $\alpha = \alpha_*$. The next result deals with the case where $y = u/\alpha_*$ and $y = F(u)$ are tangent at 0 and have no other intersection points; see Figure 1(b).

THEOREM 2.3.  *Suppose $F$ is twice differentiable at 0 and $F'(0)u > F(u)$ for $u > 0$. Then $\alpha_* = 1/F'(0)$. Suppose $I = \{i > 1 : F^{(i)}(0) \neq 0\} \neq \varnothing$. Let $\ell = \min I$. If $\alpha = \alpha_*$, then*

$$\lim_{n \to \infty} \frac{\log R_n}{\log n} = \nu_0 := \frac{2\ell - 2}{2\ell - 1}, \qquad a.s. \tag{2.10}$$

*The upper bound of $R_n$ can be strengthened to*

$$\limsup_{n \to \infty} \frac{R_n}{n^{\nu_0}(\log n)^{1-\nu_0}} \leq \left\{\frac{\ell! \sqrt{2(1+\nu_0)} F'(0)^\ell}{|F^{(\ell)}(0)|}\right\}^{2(1-\nu_0)}, \qquad a.s. \tag{2.11}$$

*Furthermore, a.s., $V_n/R_n \to \beta_*$ and true discoveries $R_n - V_n \sim (1-\beta_*)R_n$.*

The last result in Theorem 2.3 shows that choosing $\alpha = \alpha_*$ is optimal in the sense that as $n \to \infty$, the pFDR asymptotically obtains the lower bound $\beta_*$ and at the same time the number of true discoveries is unbounded, although growing sublinearly.

**3. Examples.** This section considers examples where the criticality phenomenon occurs. In each example, when a null is true, the corresponding test statistic $X \sim \Psi_0$ with density $\psi_0$; otherwise, $X \sim \Psi_1$ with density $\psi_1$. Under the random effects model, $X \sim \Psi = (1-\pi)\Psi_0 + \pi\Psi_1$. We shall focus on the right-tail $p$-value $1 - \Psi_0(X)$. From probability theory, $X \sim \Psi^{-1}(U)$, where $U \sim \text{Uniform}(0,1)$. Therefore, the $p$-value has the same distribution as $1 - \Psi_0(\Psi^{-1}(U))$, which has distribution function

$$F(u) = 1 - \Psi(\Psi_0^{-1}(1-u)) = (1-\pi)u + \pi[1 - \Psi_1(\Psi_0^{-1}(1-u))].$$

Comparing with (2.1), we get $G(u) = 1 - \Psi_1 \circ \Psi_0^{-1}(1-u)$. The density of the $p$-value is

$$(3.1) \quad F'(u) = 1 - \pi + \pi \rho(x), \qquad \text{with } x = \Psi_0^{-1}(1-u), \rho(x) = \frac{\psi_1(x)}{\psi_0(x)}.$$

From Section 2, a necessary and sufficient condition for the criticality phenomenon to occur is $\alpha_* > 0$. It can be seen that this is equivalent to $F'(0) < \infty$ (cf. Figure 1). By (3.1), $F'(0) = 1 - \pi + \pi \lim_{x \to \infty} \rho(x)$. Therefore,

$$(3.2) \qquad \text{Criticality occurs} \quad \iff \quad \lim_{x \to \infty} \rho(x) < \infty.$$

Since $\rho(x)$ is the likelihood ratio associated with $X = x$, then

$$(3.3) \quad \text{Likelihood ratio of } X \text{ is bounded} \quad \implies \quad \text{Criticality occurs}.$$

Finally, because $\Psi_0^{-1}(1-u)$ is decreasing in $u$, by (3.1) and (2.3),

$$(3.4) \quad \begin{aligned} F \text{ is strictly concave} &\iff \rho(x) \text{ is strictly increasing} \\ &\implies \alpha_* = \frac{1}{1 - \pi + \pi \lim_{x \to \infty} \rho(x)}. \end{aligned}$$

3.1. *Multiple $t$-tests.* Consider normal distributions $N(\mu_i, \sigma_i)$, $i \geq 1$, with $\mu_i$ and $\sigma_i$ being unknown. Suppose that each null is $H_i : \mu_i = 0$, and when it is false, $\mu_i = c > 0$. Assume that all $\sigma_i$ are equal but this information is unknown to the investigator. To test $H_i$, let $Y_{i,1}, \ldots, Y_{i,\nu+1}$ i.i.d. $\sim N(\mu_i, \sigma_i)$ be collected. If $\mu_i = 0$, then the $t$-statistic of the sample follows the $t$-distribution with $\nu$ degrees of freedom (d.f.); whereas if $\mu_i = c$, it follows the noncentral $t$-distribution with $\nu$ d.f. and noncentrality parameter $\delta = \sqrt{\nu+1}c/\sigma$, the density of which is

$$t_{\nu,\delta}(x) = \frac{\nu^{\nu/2}}{\sqrt{\pi}\Gamma(\nu/2)} \frac{e^{-\delta^2/2}}{(\nu+x^2)^{(\nu+1)/2}} \sum_{k=0}^{\infty} \Gamma\left(\frac{\nu+k+1}{2}\right) \frac{(\delta x)^k}{k!} \left(\frac{2}{\nu+x^2}\right)^{k/2}.$$

Observe $t_{\nu,0}(x) = t_\nu(x)$, the density of the $t$-distribution with $\nu$ d.f. We have

$$\rho(x) := \frac{t_{\nu,\delta}(x)}{t_\nu(x)}$$
(3.5)
$$= e^{-\delta^2/2} \sum_{k=0}^{\infty} \Gamma\left(\frac{\nu+k+1}{2}\right) \frac{(\delta x)^k}{k!} \left(\frac{2}{\nu+x^2}\right)^{k/2} \Big/ \Gamma\left(\frac{\nu+1}{2}\right).$$

The criticality phenomenon for the $t$-tests follows from the next result and (3.2).

PROPOSITION 3.1. *If $\delta > 0$, then*

(3.6) $$\lim_{x \to \infty} \rho(x) = e^{-\delta^2/2} \sum_{k=0}^{\infty} \Gamma\left(\frac{\nu+k+1}{2}\right) \frac{(\sqrt{2}\delta)^k}{k!} \Big/ \Gamma\left(\frac{\nu+1}{2}\right) < \infty.$$

*Furthermore, the distribution function of the p-values is strictly concave.*

PROOF. It is not hard to see (3.6) holds. For the second statement, by (3.4) it suffices to show that $\rho(x)$ is strictly increasing. It is not hard to see this is the case for $x \geq 0$. To finish the proof, denote $a = \sqrt{2}\delta$ and let

$$h_\nu(s) = \sum_{k=0}^{\infty} \Gamma\left(\frac{\nu+k+1}{2}\right) \frac{(-a)^k s^k}{k!}, \qquad \text{for } s \in [0,1].$$

Then for $x < 0$, $\rho(x) = h_\nu(\phi(x))/C$, with

$$C = e^{\delta^2/2} \Gamma\left(\frac{\nu+1}{2}\right), \qquad \phi(x) = \frac{|x|}{\sqrt{\nu+x^2}}.$$

Because $\phi$ is strictly decreasing for $x < 0$, in order to show that $\rho(x)$ is strictly increasing on $(-\infty, 0)$, it suffices to show that $h_\nu(s)$ is strictly decreasing on $(0,1)$. Since $\phi: (-\infty, 0) \to (0,1)$ is one-to-one and onto, for any $s \in (0,1)$, $h_\nu(s) = C\rho(\phi^{-1}(s)) > 0$. It is easy to see that $h'_\nu(s) = -a h_{\nu+1}(s)$. By the same argument for $h_\nu$, $h_{\nu+1}(s) > 0$ for all $s \in (0,1)$. Then $h'_\nu(s) < 0$ and hence, $h_\nu$ is strictly decreasing on $(0,1)$. □

3.2. *Multiple F-tests.* Consider regression models $Y = \boldsymbol{\beta}_i^T \mathbf{X} + \epsilon_i$, $i \geq 1$, where $\boldsymbol{\beta}_i$ is $p$-dimensional and $\epsilon_i \sim N(0, \sigma_i)$ is independent of $\mathbf{X}$. Suppose that for each $i$, the null is $H_i : \boldsymbol{\beta}_i = 0$, and when it is false, $\boldsymbol{\beta}_i = \boldsymbol{\beta} \neq 0$, with $\boldsymbol{\beta}$ being unknown. Assume that all $\sigma_i = 1$, but this information is unknown to the investigator. To test $H_i$, let an independent sample $(Y_{i,1}, \mathbf{x}_1)$, ..., $(Y_{i,\nu+p}, \mathbf{x}_{\nu+p})$ be collected, where $\mathbf{x}_1, \ldots, \mathbf{x}_{\nu+p}$ are fixed covariates for all $i$. If $\boldsymbol{\beta}_i = 0$, then the $F$-statistic of the sample follows the $F$-distribution with $(p, \nu)$ d.f.; otherwise, the statistic follows the noncentral $F$-distribution with $(p,\nu)$ d.f. and noncentrality parameter $\delta = (\boldsymbol{\beta}^T \mathbf{x}_1)^2 + \cdots + (\boldsymbol{\beta}^T \mathbf{x}_{\nu+p})^2$, the density of which is

$$f_{p,\nu,\delta}(x) = e^{-\delta/2}\theta^{p/2}x^{\nu/2-1}(1+\theta x)^{(p+\nu)/2} \sum_{k=0}^{\infty} \frac{(\delta/2)^k}{k! B(p/2+k,\nu/2)} \left(\frac{\theta x}{1+\theta x}\right)^k,$$

$$x \geq 0,$$

where $\theta = p/\nu$ and $B(a,b) = \Gamma(a)\Gamma(b)/\Gamma(a+b)$ is the Beta function. Observe that $f_{p,\nu,0}(x) = f_{p,\nu}(x)$, the density of the usual $F$-distribution with $(p,\nu)$ d.f. For $x \geq 0$,

$$(3.7) \quad \rho(x) := \frac{f_{p,\nu,\delta}(x)}{f_{p,\nu}(x)} = e^{-\delta/2} B\left(\frac{p}{2}, \frac{\nu}{2}\right) \sum_{k=0}^{\infty} \frac{(\delta/2)^k}{k! B(p/2+k,\nu/2)} \left(\frac{\theta x}{1+\theta x}\right)^k,$$

which is strictly increasing, and

$$(3.8) \quad \lim_{x\to\infty} \rho(x) = e^{-\delta/2} B\left(\frac{p}{2}, \frac{\nu}{2}\right) \sum_{k=0}^{\infty} \frac{(\delta/2)^k}{k! B(p/2+k,\nu/2)} < \infty.$$

Table 1

|  | **Monotonicity** | $\lim_{x\to\infty} \rho_k(x)$ |
|---|---|---|
| $\sigma_k = 1, \mu_k > 0$ | strictly increasing | $\infty$ |
| $\sigma_k = 1, \mu_k < 0$ | strictly decreasing | 0 |
| $\sigma_k = 1, \mu_k = 0$ | 1 | 1 |
| $\sigma_k < 1$ | maximized at $\frac{\mu_k}{1-\sigma_k^2}$ | 0 |
| $\sigma_k > 1$ | minimized at $\frac{\mu_k}{1-\sigma_k^2}$ | $\infty$ |

Therefore, the criticality phenomenon occurs with the $F$-tests. Note that by Stirling's formula, $B(a+k,b) \sim \Gamma(b)k^{-b}$ as $k \to \infty$, hence the convergence in (3.8).

3.3. *Multiple z-tests.* Suppose that the distribution under a true null is $\Psi_0 = N(0,1)$ and we have complete knowledge about this. The distribution $\Psi_1$ under a false null on the other hand is a mixture of Gaussians $p_1 N(\mu_1, \sigma_1) + \cdots + p_m N(\mu_m, \sigma_m)$, with $p_1 + \cdots + p_m = 1$ and $0 < p_k < 1$. The likelihood ratio function in this case is

$$\rho(x) = \frac{\psi_1(x)}{\psi_0(x)} = \sum_{k=1}^{m} \frac{p_k}{\sigma_k} \rho_k(x) = \sum_{k=1}^{m} \frac{p_k}{\sigma_k} \exp\left\{-\frac{(x-\mu_k)^2 - \sigma_k^2 x^2}{2\sigma_k^2}\right\}.$$

All the possibilities for each $\rho_k(x)$ are given in Table 1.

Therefore, if all the components of the Gaussian mixture have variances less than 1, from (3.2), no matter what the signs of their means are, the criticality phenomenon occurs.

3.4. *Tests involving mixtures of shift-scaled densities.* The example in Section 3.3 is a special case where the distribution under the false null is a mixture of shift-scaled versions of the distribution under the true null:

$$(3.9) \quad \Psi_1(x) = \sum_{k=1}^{m} p_k \Psi_0(s_k x - t_k), \qquad \psi_1(x) = \sum_{k=1}^{m} p_k s_k \psi_0(s_k x - t_k).$$

In most of the practical cases, the null density $\psi_0$ satisfies the tail condition

$$(3.10) \qquad \limsup_{x \to \infty} \frac{\psi_0(sx - t)}{\psi_0(x)} < \infty, \qquad \text{any } s > 1 \text{ and } t.$$

If the right tail of $\psi_0$ is rapidly decaying in the sense that

$$\lim_{x \to \infty} \frac{\psi_0(sx - t)}{\psi_0(x)} = 0, \qquad \text{any } s > 1 \text{ and } t,$$

then the criticality phenomenon occurs when $s_k < 1$ for all the components in the mixture (3.9). This is the case for Gaussian mixtures. On the other hand, if the right tail of $\psi_0$ is slowly decaying in the sense that

$$\liminf_{x \to \infty} \frac{\psi_0(sx - t)}{\psi_0(x)} = a(s, t) > 0, \qquad \text{any } s \geq 1 \text{ and } t,$$

then it is seen that (3.10) holds for any $s > 0$ and $t$. As a result, $\limsup_{x \to \infty} \rho(x) < \infty$ and the criticality phenomenon occurs.

As an example, the density of the Cauchy distribution with scale $s$ is $\mu_s(x) = \frac{1}{s\omega[1+(x/s)^2]}$, where $\omega$ denotes the circumference-diameter ratio of a circle. Suppose $X \sim \mu_1$ under a true null and $X \sim \mu_s$ with $s \neq 1$ under a false null. Because $\mu_s(x) = s^{-1}\mu_1(x/s)$ and the tail of $\mu_1(x)$ is slowly decaying, the criticality phenomenon occurs with the BH procedure when it is applied to multiple tests on the scales of the distributions.

**4. A multi-reference point procedure.** Under the random effects model, it is shown in [4] that the infimum of the pFDR attainable by multiple testing is $\underline{\beta} = (1-\pi)/\sup F'$. In general, $\sup F' < \sup \frac{F(u)}{u}$ and hence, $\underline{\beta} < \beta_*$, the infimum of the pFDR attainable by the BH procedure. This raises the question as to how to modify the BH procedure in order to improve its pFDR control.

To get some idea, consider the example of a Gaussian mixture in Section 3.3. If every component of $\Psi_1$ has a smaller variance than $\Psi_0$, then most of the $p$-values from false nulls cannot be very small. On the other hand, if the target FDR control level is small enough, the BH procedure will only reject small $p$-values, and hence, overlook most of the $p$-values from false nulls, causing a loss in power and capability to control the pFDR.

The example suggests that in order to improve power and pFDR control, an FDR controlling procedure should look at not only small $p$-values, but also moderate or even large ones for candidates of rejection. Is this reasonable? To answer the question, consider how the $p$-value is defined. First, one has to choose what amounts to "unusualness" of a single observation, for example, a large difference from 0. Then the $p$-value is defined in terms of this pre-specified unusualness under a true null. There is no guarantee that the $p$-value accounts for what is actually unusual about a population of observations from false nulls. For the Gaussian mixture example, $\Psi_1$ is unusual exactly because it generates too many "usual" observations. Rejecting small $p$-values overlooks this. However, unlike single hypothesis testing in multiple hypothesis testing the unusualness of a population of $p$-values may be detectable. In the above example, the histogram of $p$-values exhibits several peaks due to different components of $\Psi_1$. If the "attention" of an FDR controlling procedure can be distributed to the peaks, then it may have improved power and pFDR control.

4.1. *Description.* The proposed procedure will be referred to as "Procedure M," as it combines BH-type procedures at multiple locations or "reference points." Given $t \in [0,1]$, denote $R_n(t) = \#\{j \leq n : p_j \leq t\}$ and $R_n^o(t) = R_n(t-) = \#\{j \leq n : p_j < t\}$. Then the "forward" BH-type procedure at $t$ rejects nulls with $p$-values in $[t, p_{n:U_n(t)}]$, while the "backward" BH-type procedure rejects nulls with $p$-values in $[p_{n:L_n(t)}, t)$, where

$$L_n(t) = R_n^o(t) - r_n^-(t) + 1, \qquad U_n(t) = R_n(t) + r_n^+(t),$$

with

$$r_n^+(t) = \max\left\{0 \leq j \leq n - R_n(t) : p_{n:R_n(t)+j} - t \leq \frac{\alpha j}{n}\right\},$$

$$r_n^-(t) = \max\left\{0 \leq j \leq R_n^o(t) : t - p_{n:R_n^o(t)-j+1} \leq \frac{\alpha j}{n}\right\}.$$

The total number of rejections at $t$ is then $r_n(t) = r_n^+(t) + r_n^-(t) = U_n(t) - L_n(t) + 1$.

Given reference points $0 = t_0 < t_1 < \cdots < t_M = 1$, Procedure M is as follows:

STEP 1. Select $t_{i_1} < t_{i_1} < \cdots < t_{i_m}$, such that:

(1) $r_n(t_{i_k}) \geq (\log n)^c$ for each $k = 1, \ldots, m$, where $c > 1$ is a parameter;
(2) $[L_n(t_{i_k}), U_n(t_{i_k})]$ are disjoint, that is, $U_n(t_{i_k}) < L_n(t_{i_{k+1}})$ for $k = 1, \ldots, m$; and
(3) $\sum_{k=0}^m r_n(t_{i_k})$ is the largest among all subsets of $t$'s satisfying (1) and (2).

STEP 2. Reject all $p_j$ falling into the union of $J(t_{i_k}) := [p_{n:L_n(t_{i_k})}, p_{n:U_n(t_{i_k})}]$.

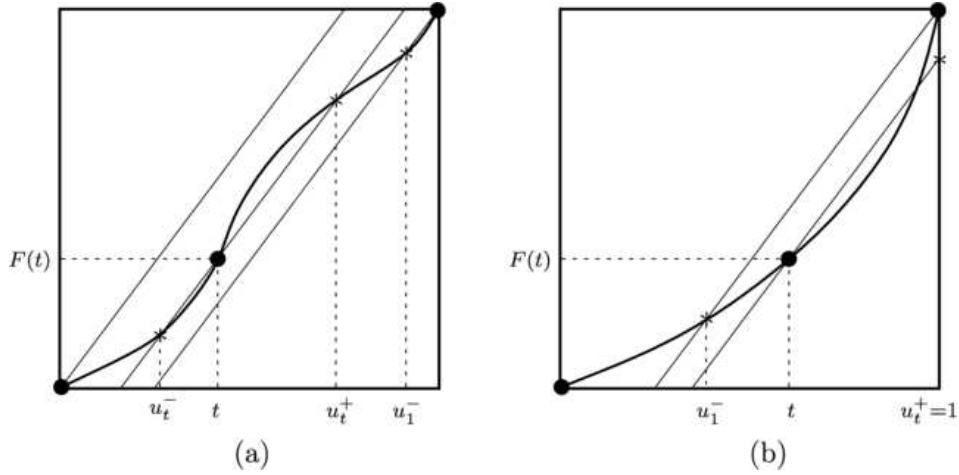

Fig. 2. *Rejections at different reference points. The parallel thin lines passing through $(0,0)$, $(t, F(t))$ and $(1,1)$ all have slope $1/\alpha$.*

Condition (1) in Step 1 requires that every selected reference point should have at least $(\log n)^c$ rejections. As for the BH procedure, the BH-type procedures at each $t$ may have a critical value $\alpha_*(t) > 0$ such that, for $\alpha < \alpha_*(t)$, $r_n(t)/(\log n)^c \xrightarrow{\text{a.s.}} 0$ for any $c > 1$ and for $\alpha \geq \alpha_*(t)$, $r_n(t)$ is of order $n^a$ for some $a > 0$. Condition (1) is used in the case $\alpha < \alpha_*(t)$.

The optimization in Step 1 can be efficiently computed by dynamic programming; see the Appendix. In general, Procedure M generates several random rejection intervals, each one a connected component of the union of $J(t_{i_k})$ associated with selected $t_{i_k}$. The BH procedure, on the other hand, always generates a single random rejection interval.

4.2. *Some justifications.* To see roughly why Procedure M can increase power and improve pFDR control, first consider Figure 2(a). In this case, by Theorem 2.1, the BH procedure asymptotically has power 0. On the other hand, for the $t$ in Figure 2(a), $\alpha > 1/F'(t)$. By Proposition 4.1,

$$
\begin{aligned}
p_{n:L_n(t)} &\xrightarrow{\text{a.s.}} u_t^- := \inf\{x \in [0,t] : t - x \leq \alpha(F(t) - F(x))\}, \\
p_{n:U_n(t)} &\xrightarrow{\text{a.s.}} u_t^+ := \sup\{x \in [t,1] : x - t \leq \alpha(F(x) - F(t))\}.
\end{aligned}
\tag{4.1}
$$

Similar to the BH procedure, among the nulls with $p$-values between $p_{n:L_n(t)}$ and $p_{n:U_n(t)}$, the expected proportion of true nulls is less than $\alpha$. Therefore, the BH-type procedure at $t$ has positive power. Likewise, the backward BH-type procedure at $t = 1$ has positive power. Procedure M thus controls the FDR at $\alpha$ with positive power.

Now imagine that $\tilde{p}_j = 1 - p_j$ instead of $p_j$ are used for the tests. Then the distribution function of the $p$-values becomes $\tilde{F}(x) = 1 - F(1-x)$. Since $\tilde{F}'(0) = F'(1) > 1/\alpha$, by Theorem 2.2 the BH procedure has positive power. However, by combining the BH-type procedure at $1-t$, the power is increased by a positive factor.

To see why Procedure M requires selection of $t_{i_j}$, note that $[p_{n:L_n(t)}, p_{n:U_n(t)}]$ may overlap for different $t$. In Figure 2(b), since $F$ is convex, for $t \in (u_1^-, 1)$, $u_t^- = t$ and $u_t^+ = 1$, whereas for $t < u_1^-$, $u_t^- = u_t^+ = t$. Therefore, for $n \gg 1$, only $t_M = 1$ is selected.

As an example, consider the tests on Cauchy distributions in Section 3.4 again. If the $p$-value of each observation $X$ is defined as the two-tailed probability $\Pr\{|W| > |X|\}$, where $W \sim \mu_1$ is independent of $X$, then it has distribution function

$$F_{\text{2-tail}}(u) = (1-\pi)u + \frac{2\pi}{\omega}\cot^{-1}\left(\frac{1}{s}\cot\frac{\omega u}{2}\right),$$

with $F'_{\text{2-tail}}(0) = 1 - \pi + \pi s$.

The critical value is $\alpha_* = \frac{1}{1-\pi+\pi s}$. Let $\alpha \in (\alpha_*, 1)$. As $F_{\text{2-tail}}$ is strictly concave, the BH procedure and Procedure M are asymptotically the same, both rejecting $p$-values in $[0, u_*]$, with $u_* \in (0, 1)$ the solution to $u/\alpha = F_{\text{2-tail}}(u)$. On the other hand, if the $p$-value is defined as the left-tail probability $\Pr\{W < X\}$, then its distribution function is

$$F_{\text{left}}(u) = \begin{cases} \frac{1}{2}F_{\text{2-tail}}(2u), & 0 \le u \le \frac{1}{2}, \\ 1 - \frac{1}{2}F_{\text{2-tail}}(2-2u), & \frac{1}{2} \le u \le 1. \end{cases}$$

Asymptotically, the BH procedure only rejects $p$-values in $[0, u_0^+]$, while Procedure M rejects those in $[0, u_0^+] \cup [1 - u_1^-, 1]$, with $u_0^+ \in (0, 1)$ the solution to $u/\alpha = F_{\text{left}}(u)$ and $u_1^- \in (0, 1)$ the solution to $(1-u)/\alpha = 1 - F_{\text{left}}(u)$. It is seen that $u_0^+ = 1 - u_1^- = u_*/2$. Thus, the power of the BH procedure is reduced by half, while that of Procedure M is unchanged. Procedure M therefore is less sensitive to the choice of $p$-values.

4.3. *Theoretical properties.* The following characterizations of $L_n$ and $U_n$ by continuous stopping times and fixed points generalize those in [17] and [8].

PROPOSITION 4.1. $L_n(t) = R_n^o(T_n^-(t)) + 1$ and $U_n(t) = R_n(T_n^+(t))$, where

(4.2)
$$T_n^-(t) = \inf\left\{x \le t : \frac{t-x}{\alpha} \le \frac{[R_n^o(t) - R_n^o(x)] \vee 1}{n}\right\},$$
$$T_n^+(t) = \sup\left\{x \ge t : \frac{x-t}{\alpha} \le \frac{[R_n(x) - R_n(t)] \vee 1}{n}\right\}.$$

*Therefore, Step 2 of Procedure M is the same as rejecting all nulls with $p$-values in $[T_n^-(t_{i_k}), T_n^+(t_{i_k})]$, $k = 0, \ldots, m$. Furthermore, $T_n^\pm(t) \xrightarrow{\text{a.s.}} u_t^\pm$ and (4.1) holds.*

The results on power and pFDR control in Section 2 for the BH procedure can be generalized to Procedure M. Let $p_*^+(t) = F(u_t^+) - F(t)$, $p_*^-(t) = F(t) - F(u_t^-)$.

THEOREM 4.1. *Suppose $t_i - \alpha F(t_i)$ are different from each other:*

(1) *If $p_*^\pm(t_i) = 0$ and $F'(t_i) < 1/\alpha$ for all $i$, then a.s., for $n \gg 1$, $R_n = 0$. Furthermore,*

$$(r_n^+(0), r_n^-(t_1), r_n^+(t_1), \ldots, r_n^-(t_{M-1}), r_n^+(t_{M-1}), r_n^-(1))$$
$$\xrightarrow{d} (\tau_0, \tilde\tau_1, \tau_1, \ldots, \tilde\tau_{M-1}, \tau_{M-1}, \tilde\tau_M), \qquad as\ n \to \infty,$$

*where the $\tau$'s and $\tilde\tau$'s are independent, each $\tau_k$ and $\tilde\tau_k$ following the distribution of the last excursion time into $(-\infty, 0)$ of the random walk $S_0 = 0$, $S_k = S_{k-1} + \gamma_k - \alpha F'(t_k)$, with $\gamma_1, \gamma_2, \ldots$ i.i.d. $\sim \text{Exp}(1)$.*

(2) *If $p_*^+(t_i) + p_*^-(t_i) > 0$ for some $i = 0, \ldots, M$, then a.s.,*

$$\limsup_{n\to\infty} \text{FDR} = \limsup_{n\to\infty} \text{pFDR} \le (1-\pi)\alpha,$$

$$\lim_{n\to\infty} \text{power}_n = \left(\frac{1-\alpha}{\pi} + \alpha\right)\Pi,$$

*where*

$$\Pi = \max_{\substack{S \subset \{t_0, \ldots, t_M\}: [u_s^-, u_s^+] \\ \text{are disjoint for } s \in S}} \left\{\sum_{s \in S} [p_*^-(s) + p_*^+(s)]\right\}.$$

(3) *Suppose for each $i$, $I_i := \{k > 1 : F^{(k)}(t_i) \ne 0\} \ne \varnothing$. Let $\ell_i = \min I_i$. If $p_*^\pm(t_i) = 0$ for all $i$ but $F'(t_i) = 1/\alpha$ for at least one of them, then*

$$\lim_{n\to\infty} \text{FDR} = \lim_{n\to\infty} \text{pFDR} = (1-\pi)\alpha, \qquad \frac{\log R_n}{\log n} \xrightarrow{a.s.} \frac{2\ell - 2}{2\ell - 1},$$

*where $\ell = \max\{\ell_i : F'(t_i) = 1/\alpha\}$. Additionally, a.s., for $n \gg 1$, the set of rejected p-values consists exactly of those in $[T_n^-(t_i), T_n^+(t_i)]$ with $F'(t_i) = 1/\alpha$.*

Note that, unlike the BH procedure, if Procedure M cannot control the pFDR, then for $n \gg 1$ it will make no rejections. This way of signaling the controllability of the pFDR can also be used to control other types of error rates for multiple testing [4].

**5. Numerical study.** In this section we report simulation studies on the criticality phenomenon of the BH procedure and Procedure M. All the simulations are conducted using R [11]. In the simulations, $p$-values are sampled as follows:

Table 2

(A) *simulations for the BH procedure with* $\alpha \equiv 0.3 < \alpha_*$. $t_{\nu,\delta}$ *is the noncentral t-distribution with* $\nu$ *d.f. and parameter* $\delta$. $F_{m,n,\delta}$ *is the noncentral F-distribution with* $(m,n)$ *d.f. and parameter* $\delta > 0$. $\widehat{\mathrm{FDR}}$ *and* $\widehat{\mathrm{pFDR}}$ *are the MC estimates of* $\mathrm{FDR} = (1-\pi)\alpha$ *and* $\mathrm{pFDR} = (1-\pi)\alpha_*$. (B) *Estimated* $\hat{p}_k = \widehat{\mathrm{Pr}}\{\tau = k\}$ *from the first simulation in* (A) *vs.* $p_k = \mathrm{Pr}\{\tau = k\}$ *based on* (2.5)

| | | | | (A) | | | | |
|---|---|---|---|---|---|---|---|---|
| Simu. | $\Psi_0$ | $\Psi_1$ | $\pi$ | $\alpha_*$ | $\widehat{\mathrm{FDR}}$ | FDR | $\widehat{\mathrm{pFDR}}$ | pFDR |
| 1 | $t_{10}$ | $t_{10,1}$ | 0.05 | 0.512 | 0.272 | 0.285 | 0.532 | 0.486 |
| 2 | $t_{10}$ | $t_{10,1}$ | 0.02 | 0.724 | 0.297 | 0.294 | 0.749 | 0.709 |
| 3 | $F_{10,10}$ | $F_{10,10,3}$ | 0.05 | 0.892 | 0.294 | 0.285 | 0.867 | 0.847 |
| 4 | $F_{10,10}$ | $F_{10,10,3}$ | 0.02 | 0.954 | 0.305 | 0.294 | 0.941 | 0.935 |

| | | | | (B) | | | | |
|---|---|---|---|---|---|---|---|---|
| $k$ | $\hat{p}_k$ | $p_k$ | $k$ | $\hat{p}_k$ | $p_k$ | $k$ | $\hat{p}_k$ | $p_k$ |
| 0 | 0.4900 | 0.4140 | 8 | 0.0176 | 0.0221 | 16 | 0.0022 | 0.0060 |
| 1 | 0.1474 | 0.1350 | 9 | 0.0098 | 0.0185 | 17 | 0.0022 | 0.0051 |
| 2 | 0.0908 | 0.0881 | 10 | 0.0094 | 0.0155 | 18 | 0.0016 | 0.0044 |
| 3 | 0.0604 | 0.0646 | 11 | 0.0066 | 0.0131 | 19 | 0.0004 | 0.0038 |
| 4 | 0.0474 | 0.0500 | 12 | 0.0068 | 0.0112 | 20 | 0.0008 | 0.0033 |
| 5 | 0.0410 | 0.0400 | 13 | 0.0042 | 0.0095 | 21 | 0.0010 | 0.0029 |
| 6 | 0.0296 | 0.0323 | 14 | 0.0038 | 0.0081 | 22 | 0.0004 | 0.0025 |
| 7 | 0.0198 | 0.0266 | 15 | 0.0038 | 0.0070 | 23 | 0.0010 | 0.0022 |

1. Sample $\theta \sim \mathrm{Bernoulli}(\pi)$, where $\pi$ is the population fraction of false nulls.
2. If $\theta = 0$, sample $p \sim \mathrm{Unif}(0,1)$ and return $p$.
3. If $\theta = 1$, sample $X \sim \Psi_1$ and return $p = 1 - \Psi_0(X)$, where $\Psi_0$ is the distribution under a true null, and $\Psi_1$ the distribution under a false null.

5.1. *Simulation study on the BH procedure with* $\alpha < \alpha_*$. This part of the study consists of four simulations on $t$- and $F$-tests. Since the distributions of $p$-values are strictly concave, we apply (3.4), (3.6) and (3.8) to compute the critical value $\alpha_*$. The parameters of the null distributions and values of $\alpha_*$ are shown in Table 2(A), columns 2–5. The Appendix has some remarks on the evaluation of $\alpha_*$.

In the simulations, $\alpha = 0.3$. Each simulation contains 5000 runs. As the distribution of $R_n$ appears to converge slowly, in each run the BH procedure is applied to $n = 10^5$ sample $p$-values. The FDR and pFDR are estimated by the Monte Carlo (MC) average

$$\widehat{\mathrm{FDR}} = \frac{1}{5000} \sum_{j=1}^{5000} \frac{v_j}{r_j \vee 1}, \qquad \widehat{\mathrm{pFDR}} = \frac{1}{N} \sum_{j=1}^{5000} \frac{v_j \mathbf{1}\{r_j > 0\}}{r_j \vee 1},$$

where $r_j$ and $v_j$ are the numbers of rejections and false rejections in the $j$th run, respectively, and $N$ the number of runs with $r_j > 0$. From the last four columns of Table 2(A), we see good agreement between the simulations and the theoretical results that $\text{pFDR} \to (1-\pi)\alpha_*$ and $\text{FDR} \equiv (1-\pi)\alpha$.

We next compare the limiting distribution of $R$ in (2.5) and its estimate

$$\hat{p}_k = \frac{1}{5000} \sum_{j=1}^{5000} \mathbf{1}\{r_j = k\}.$$

Table 2(B) collects the results of the first simulation, which has $\Psi_0 = t_{10}$, $\Psi_1 = t_{10,1}$ and $\pi = 0.05$. For relatively small values of $k$, $\hat{p}_k$ agrees with $p_k = \Pr\{\tau = k\}$ reasonably well. The Kullback–Leibler (KL) distance between $\hat{p}$ and $p$, $D(\hat{p}\|p) = \sum_{i=0}^{\infty} \hat{p}_k \log \frac{\hat{p}_k}{p_k}$, is equal to 0.0409. On the other hand, the total variation (TV) distance $d_{\text{TV}}(\hat{p}, p) = 0.1847$, which is the largest among the four simulations. The KL and TV distances from the other three simulations are $(0.0042, 0.0425)$, $(0.0018, 0.0220)$ and $(0.0020, 0.0306)$, respectively, indicating that the smaller $\alpha/\alpha_*$ is, the faster $\hat{p}$ converges to $p$.

5.2. *Simulation study on the BH procedure with $\alpha > \alpha_*$.* This part of the study consists of simulations on $t$-, $F$- and $z$-tests. The parameters of the simulations are shown in Table 3(A), columns 2–4.

In each simulation, $\alpha = 0.25$, the number of runs is 400, and each run samples $n = 40000$ $p$-values. The MC estimates of pFDR and the theoretical values of FDR are shown in columns 5–6, which agree with each other very well. In all the simulations, $\widehat{\text{FDR}}$ is identical to $\widehat{\text{pFDR}}$.

We next numerically evaluate the LIL in Theorem 2.2. In principle, this requires that in each run we sample an infinite sequence of $p$-values, apply the BH procedure to the first $n = 1, 2, \ldots$ of them, and then check $\limsup_n \pm L_n = \lambda$, where we denote

$$(5.1) \qquad L_n = \frac{R_n - np_*}{\sqrt{n \log \log n}}, \qquad \lambda = \frac{\sqrt{2p_*(1 - p_*)}}{1 - \alpha F'(u_*)}.$$

However, this is problematic because (1) the convergence is very slow, and (2) applying the BH procedure for every $n$ is computationally too costly. We instead take the following approach. In each run, sample $N = 320^2 = 102400$ $p$-values. Set $n_k = \lfloor (200 + 0.2k)^2 \rfloor$. Then, for $k = 0, 1, \ldots, 600$, apply the BH procedure to the first $n_k$ $p$-values to get $L_{n_k}$. Finally, the output of the run is $\max_k |L_{n_k}|$.

Due to the slow convergence of the LIL, we do not expect to observe convergence of $\max_k |L_{n_k}|$ to $\lambda$ in the simulations. However, if we can show that, for most of the runs, $\max_k |L_{n_k}| \leq \lambda$, then it empirically demonstrates $R_n = np_* + O_p(\sqrt{\log \log n / n})$.

The results of the simulations are summarized in Table 3(B). The values of $p_*$, $F'(u_*)$ and $\lambda$ involved in the LIL (5.1) are shown as well. Note that $p_*$ is also the asymptotic limit of the proportion of rejected $p$-values; see (2.4). The Appendix has some comments on its numerical evaluation. For $t \geq 1$, let $\hat{Q}(t)$ be the proportion of runs with $\max_k |L_{n_k}| > t\lambda$. It is seen that, for most of the runs, $\max_k |L_{n_k}| \leq \lambda$. Figure 3 plots all the curves of $L_n/\lambda$ versus $\sqrt{n \log \log}$ in the 400 runs of simulation 5. Again, it shows that, for most of the time, $L_n/\lambda$ stays in $[-1, 1]$. The plots from the other five simulations show the same property.

5.3. *Simulation study of Procedure* M. This part consists of six simulations to compare the BH procedure and Procedure M. In Section 4 it is noted that for Procedure M, it is a useful contraint for pFDR control that each selected reference point generate at least $(\log n)^c$ rejections. To examine this, we test Procedure M for $c = 1.5$, 2, and also its modified version which only requires that each selected reference point generate at least one rejection. These three procedures are denoted M(1.5), M(2) and M(0), respectively.

Table 3

(A) *Parameters in the simulations for the BH procedure with* $\alpha \equiv 0.25 > \alpha_*$. *The MC estimate* $\widehat{\text{pFDR}}$ *is obtained for the total number of p-values per run* $n = 40000$.

(B) *Simulation results on the LIL for the simulations in* (A). $\hat{Q}(t)$ *is the proportion of runs with* $\max_k |L_{n_k}| > t\lambda$

| | | | (A) | | | |
|---|---|---|---|---|---|---|
| Simu. | $\Psi_0$ | $\Psi_1$ | $\pi$ | $\alpha_*$ | $\widehat{\text{pFDR}}$ | FDR |
| 1 | $t_{20}$ | $t_{20,2}$ | 0.2 | $1.5 \times 10^{-3}$ | 0.2005 | 0.2 |
| 2 | $t_{20}$ | $t_{20,2}$ | 0.15 | $2 \times 10^{-3}$ | 0.2122 | 0.2125 |
| 3 | $F_{20,30}$ | $F_{20,30,30}$ | 0.2 | $6 \times 10^{-5}$ | 0.2001 | 0.2 |
| 4 | $F_{20,30}$ | $F_{20,30,30}$ | 0.15 | $8 \times 10^{-4}$ | 0.2122 | 0.2125 |
| 5 | $N(0,1)$ | $N(2,1)$ | 0.2 | 0 | 0.2002 | 0.2 |
| 6 | $N(0,1)$ | $N(2,1)$ | 0.15 | 0 | 0.2127 | 0.2125 |

| | | | (B) | | | |
|---|---|---|---|---|---|---|
| Simu. | $p_*$ | $F'(u_*)$ | $\lambda$ | $\hat{Q}(1)$ | $\hat{Q}(1.1)$ | $\hat{Q}(1.2)$ | $\hat{Q}(1.5)$ |
| 1 | 0.1330 | 1.970 | 0.946 | 0.105 | 0.065 | 0.043 | 0.005 |
| 2 | 0.0855 | 2.123 | 0.843 | 0.113 | 0.060 | 0.040 | 0.005 |
| 3 | 0.1876 | 1.513 | 0.888 | 0.125 | 0.068 | 0.025 | 0 |
| 4 | 0.1313 | 1.649 | 0.813 | 0.155 | 0.108 | 0.065 | 0.015 |
| 5 | 0.1453 | 1.781 | 0.898 | 0.125 | 0.073 | 0.040 | 0.008 |
| 6 | 0.0974 | 1.896 | 0.797 | 0.1 | 0.053 | 0.038 | 0.008 |

Each simulation contains 1000 runs. In each run, the BH procedure and M($c$), $c = 0$, 1.5, 2, are applied to the same set of 20000 sample $p$-values. The distribution $\Psi_0$ when a null is true is $N(0,1)$, the target FDR control level is $\alpha = 0.3$, and the reference points for M($c$) are $0.01k$, $0 \le k \le 100$. The other parameters for the simulations and the critical values $\alpha_*$ for the BH procedure are shown in Table 4(A), columns 2–4.

The results are summarized in Table 4(B). In simulation 1, because $\alpha$ is only about a third of $\alpha_*$, the BH procedure has little power and pFDR $\approx 1$. On the other hand, $\sup F' = 4.83$, where $F'$ is the density of the $p$-values; see Table 4(A). From (4.1) and Theorem 4.1, if Procedure M has reference points close to where $F'$ is maximized, then it can attain pFDR $\approx (1-\pi)\alpha = 0.285$ with positive power. This is the case for M(1.5) and M(2). In contrast, M(0) has substantially higher pFDR despite a little more power. This demonstrates that by introducing some constraint on the minimum number of rejections at each selected reference point, Procedure M can attain pFDR control with positive power even when the BH procedure cannot.

In simulation 2, since $\sup F' = 2.53$, $\alpha < 1/F'(t)$ for every $t$. From Theorem 4.1, Procedure M cannot attain pFDR $= 0.3$. M(1.5) signals this with small $P(R_n > 0) \approx \frac{\widehat{\text{FDR}}}{\widehat{\text{pFDR}}} \approx 0.03$. M(2) signals this by making no rejections. On the other hand, the BH procedure and M(0) fail to signal the fact that pFDR $= 0.3$ is not attainable.

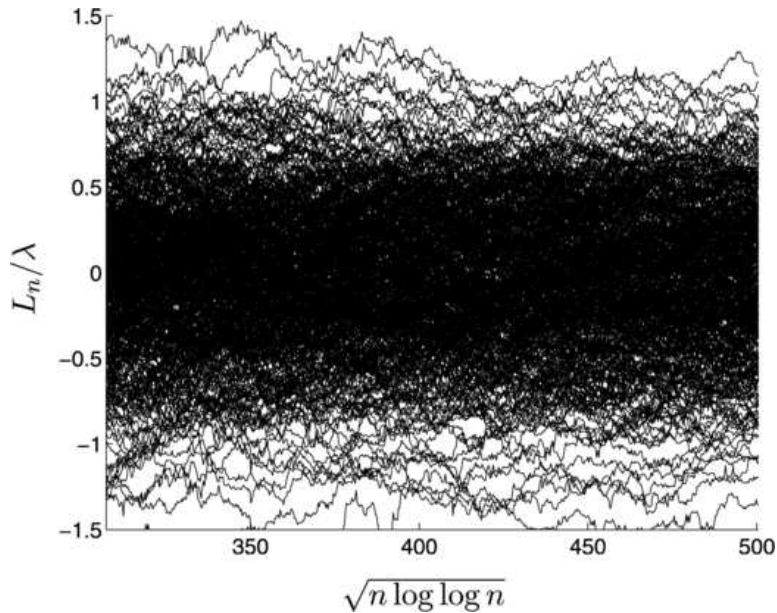

FIG. 3. *Plots of $L_n/\lambda$ vs. $\sqrt{n \log \log n}$ for all the runs in simulation 5 in Table 3(A).*

In simulations 3–6, the distribution $\Psi_1$ has a component with SD greater than that of $\Psi_0$. By Section 3.3, the BH procedure can attain any level of pFDR with positive power. However, in simulations 3–4, $u_*$ is extremely small; see Table 4(A). Recall that $u_*$ is the limit of the largest $p$-value rejected by the BH procedure; see (2.4). Consequently, with overwhelming chance, the BH procedure only rejects $p$-values from $\Psi_0$, and hence, has little power. In contrast, M(1.5) can attain pFDR $\approx (1-\pi)\alpha$ in both simulations 3 and 4, while M(2) has substantially worse performance than M(1.5) in simulation 4, indicating that its constraint on the minimum number of rejections per selected reference point is too hard when $\pi$ is small.

Finally, in simulations 5–6, $p_*$ is large enough for the BH procedure to attain pFDR control with moderate power. In simulation 5, both M(1.5) and M(2) can attain pFDR control with almost 3 times more power than the BH procedure. However, in simulation 6, while M(1.5) has similar performance to that of the BH procedure, M(2) is significantly worse.

TABLE 4
(A) *Parameters in the simulations for the BH procedure vs. Procedure* M *with* $\alpha \equiv 0.3$. $F'$ *is the density of the upper-tail $p$-value under* $\Psi_0 = N(0,1)$ *for* $X \sim (1-\pi)\Psi_0 + \pi\Psi_1$. $u_*$ *is the limit of the largest $p$-value rejected by the BH procedure.* (B) *Simulation results*

| (A) | | | | | |
|---|---|---|---|---|---|
| Simu. | $\Psi_1$ | $\pi$ | $\alpha_*$ | $\sup F'$ | $u_*$ |
| 1 | $\frac{1}{2}N(-1.3, 0.015) + \frac{1}{2}N(1, 0.015)$ | 0.05 | 0.9107 | 4.8307 | 0 |
| 2 | Same as 1 | 0.02 | 0.9623 | 2.5323 | 0 |
| 3 | $\frac{2}{5}N(-1.3, 0.015) + \frac{2}{5}N(1, 0.015) + \frac{1}{5}N(-4, 2)$ | 0.05 | 0 | $\infty$ | $1.4 \times 10^{-9}$ |
| 4 | Same as 3 | 0.02 | 0 | $\infty$ | $2.9 \times 10^{-10}$ |
| 5 | $\frac{2}{5}N(-1.3, 0.015) + \frac{2}{5}N(1, 0.015) + \frac{1}{5}N(4, 2)$ | 0.05 | 0 | $\infty$ | 0.01 |
| 6 | Same as 5 | 0.02 | 0 | $\infty$ | 0.0039 |

| (B) | | | | | | | | | |
|---|---|---|---|---|---|---|---|---|---|
| Simu. | Type | $\widehat{\text{FDR}}$ | $\widehat{\text{pFDR}}$ | $\widehat{\text{Power}}$ | Simu. | Type | $\widehat{\text{FDR}}$ | $\widehat{\text{pFDR}}$ | $\widehat{\text{Power}}$ |
| 1 | BH | 0.29 | 1 | 0 | 2 | BH | 0.296 | 1 | 0 |
|  | M(0) | 0.3541 | 0.3541 | 0.7647 |  | M(0) | 0.9559 | 0.9559 | 0.0145 |
|  | M(1.5) | 0.2882 | 0.2882 | 0.7643 |  | M(1.5) | 0.0145 | 0.4823 | 0.0015 |
|  | M(2) | 0.2882 | 0.2882 | 0.7641 |  | M(2) | 0 | NA | 0 |
| 3 | BH | 0.2907 | 0.9954 | $2.9 \times 10^{-6}$ | 4 | BH | 0.315 | 1 | 0 |
|  | M(0) | 0.3835 | 0.3835 | 0.5498 |  | M(0) | 0.7103 | 0.7103 | 0.1459 |
|  | M(1.5) | 0.2924 | 0.2924 | 0.5471 |  | M(1.5) | 0.2912 | 0.2913 | 0.1365 |
|  | M(2) | 0.2917 | 0.2917 | 0.5356 |  | M(2) | 0.0242 | 0.3565 | 0.0112 |
| 5 | BH | 0.2821 | 0.2821 | 0.1475 | 6 | BH | 0.2911 | 0.2911 | 0.1360 |
|  | M(0) | 0.3835 | 0.3835 | 0.5432 |  | M(0) | 0.7100 | 0.7100 | 0.1454 |
|  | M(1.5) | 0.2923 | 0.2923 | 0.5404 |  | M(1.5) | 0.2913 | 0.2913 | 0.1362 |
|  | M(2) | 0.2914 | 0.2914 | 0.5283 |  | M(2) | 0.0217 | 0.3559 | 0.0101 |

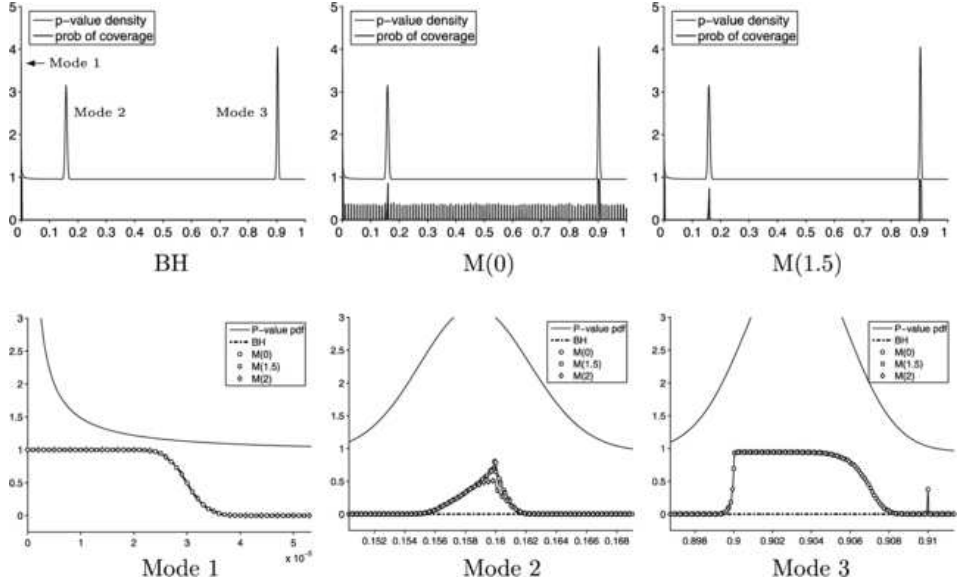

Fig. 4. *Simulation 5 for the BH procedure and* M$(c)$, $c = 0, 1.5, 2$. *Top:* $F'(t)$ *vs.* $P_{\text{cover}}(t)$. *Bottom: enlarged views of the functions around the three modes.*

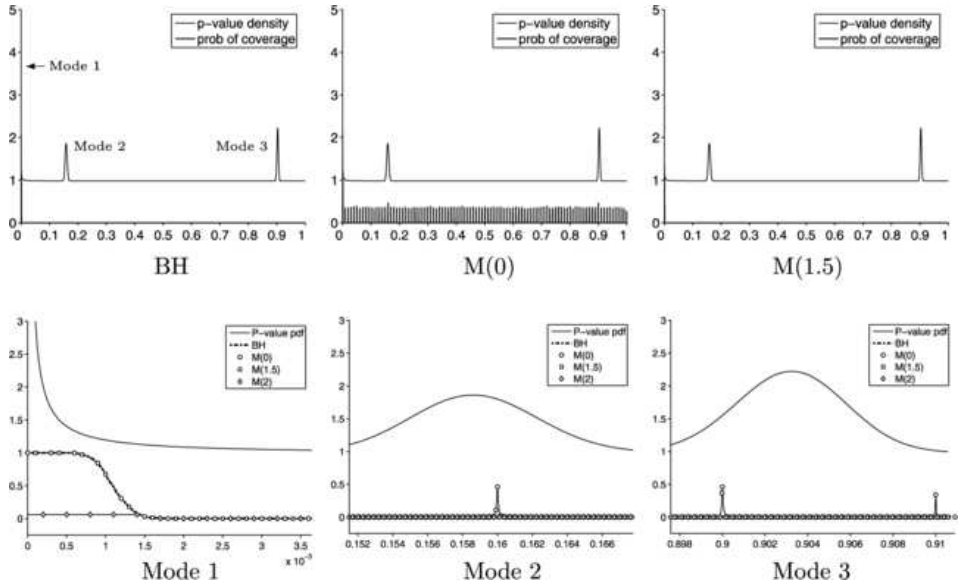

Fig. 5. *Simulation 6 for the BH procedure and* M$(c)$, $c = 0, 1.5, 2$. *Top:* $F'(t)$ *vs.* $P_{\text{cover}}(t)$. *Bottom: enlarged views of the functions around the three modes.*

To identify the sources of the differences among the procedures, one way is to look at what $p$-values are likely to be rejected. We consider $P_{\text{cover}}(t)$, the probability that $t \in [0, 1]$ is covered by a rejection region, which can be estimated by

$$P_{\text{cover}}(t) = \begin{cases} \frac{1}{1000} \sum_{i=1}^{1000} \mathbf{1}\{t \leq \text{largest rejected } p\text{-value}\}, & \text{for BH}, \\ \frac{1}{1000} \sum_{i=1}^{1000} \mathbf{1}\{t \in J(s) \text{ for some selected reference point } s\}, & \\ & \text{for M}(c), \end{cases}$$

where $J(s)$ is defined as in Procedure M, Step 2.

Figures 4 and 5 plot $P_{\text{cover}}(t)$ as well as the density of $p$-values $F'(t)$ for simulations 5 and 6, respectively. The density has three modes in each case. In simulation 5, for the BH procedure, $P_{\text{cover}}(t)$ has a single peak aligned with the mode at 0, and for M(1.5), $P_{\text{cover}}(t)$ has three peaks aligned with the three modes, which explains why it has significantly more power than the BH procedure. The shape of $P_{\text{cover}}(t)$ for M(2) is very close to that for M(1.5). In contrast, for M(0), $P_{\text{cover}}(t)$ has many peaks located at the reference points in addition to the three major peaks. As a result, M(0) rejects too many $p$-values that are highly likely to be associated with true nulls, which substantially decreases its capability to control the pFDR. For all the procedures, $P_{\text{cover}}(t)$ is almost identical around 0.

The situation is different in simulation 6. Because the values of $F'$ at the nonzero modes are less than 3, neither M(1.5) nor M(2) makes rejections around them. Around 0, the BH procedure, M(0) and M(1.5) yield almost identical $P_{\text{cover}}(t)$, whereas M(2) yields one well below the others. The BH procedure and M(1.5) therefore have similar power and pFDR control, and notably outperform the other two.

**6. Discussion.** Exploratory studies aim to identify real novel signals from a large number of observations. For any procedure used to detect the signals, its performance is constrained by the amount of useful information in the data. The BH procedure is no exception. In this article we demonstrate a criticality phenomenon that imposes constraints on the procedure's power and pFDR control. Roughly speaking, the criticality is not only due to the bounded likelihood ratios, or "signal-noise ratios" of the data [see equation (3.3)], but also may be due to the procedure missing important clues that separate false nulls from true ones. The whereabouts of these clues, by the nature of exploratory study, is unknown a priori. Based on this perspective, we propose a multi-reference point procedure which conducts testing at multiple locations in the domain of $p$-values in order to catch useful clues.

Many questions remain to be answered on what constraints multiple testing may have and how to tackle them. First, although the random effects model adapted here is helpful in identifying some of the constraints, how dependence among observations may affect power and pFDR control is yet to be seen. Second, how to detect the performance bound for a procedure. From Section 2, to estimate a lower bound for the pFDR attainable by the BH procedure, one possible way is to utilize an estimated distribution function of the $p$-values. Alternatively, different subsamples of the $p$-values may be tested at a given target control level. The distribution of the number of rejections then could give some clue as to whether the target control level is below the infimum of the attainable pFDR.

The multi-reference point procedure has much room for improvement. In essence, the procedure itself consists of multiple tests, one per reference point. Therefore, some regulations are needed for the procedure, otherwise it may have the same types of problems it is intended to address. One issue is how to determine the minimum number of rejections for each selected reference point. Second, when the number of tested hypotheses increases, it is reasonable to increase the number of reference points. This raises the question as to whether there is an optimal rate of increase. It is also possible that the reference points can be better allocated according to an estimated density function of the $p$-values.

## APPENDIX: SOME NUMERICAL ISSUES

*Procedure* M. The optimization in Step 1 is computed by dynamic programming. First, remove all $t_j$ with $r(t_j) < (\log n)^c$. Then, relabel the remaining reference points as $s_1, \ldots, s_k$ so that $U_n(s_1) \leq U_n(s_2) \leq \cdots \leq U_n(s_k)$. Denote $l_j = L_n(s_j)$, $u_j = U_n(s_j)$. Let $I_j = \{l_j, l_j + 1, \ldots, u_j\}$. Step 1 of Procedure M requires a solution to

$$S = \max\left\{\sum_{a \in A} |I_a| : A \subset \{1, \ldots, k\} \text{ such that } I_a, a \in A, \text{ are disjoint}\right\},$$

as well as the maximizing $A$. Let $S_0 = 0$ and for $j = 1, \ldots, k$,

$$S_j = \max\left\{\sum_{a \in A} |I_a| : A \subset \{1, \ldots, j\} \text{ such that } I_a, a \in A, \text{ are disjoint}\right\}.$$

Then $S = S_k$ and $S_j = \max\{S_{i(j)} + |I_j|, S_{j-1}\}$, where $i(j) = \max\{s : u_s < l_j\}$. This recursion leads to the following dynamic programming procedure:

for $j = 1, \ldots, k$
   if $S_{i(j)} + |I_j| > S_{j-1}$ then: select$[j] \leftarrow 1$, previous$[j] \leftarrow i(j)$, $S_j \leftarrow S_{i(j)} + |I_j|$
   else: select$[j] \leftarrow 0$, previous$[j] \leftarrow j - 1$, $S_j \leftarrow S_{j-1}$
$A \leftarrow \varnothing$, $s \leftarrow k$

```
    while s > 0
        if select[s] = 1 then: A ← A ∪ {s}, s ← previous[s]
        else: s ← s − 1
return  A and S = S_k
```

*Numerical evaluation of $\alpha_*$ and $p_*$.* For $t$-tests, $\alpha_*$ is evaluated via (3.4) and (3.6). To improve numerical precision, each term in (3.6) can be evaluated by $\exp(z_k)$ with $z_k = g((\nu+k_1)/2) - g(k+1) - g((\nu+1)/2) - k\log(\sqrt{2}\delta)$, where $g(x) = \log\Gamma(x)$ is evaluated by `lgamma` in R. For $F$-tests, $\alpha_*$ is evaluated via (3.4) and (3.8). Each term in (3.8) can be computed by $\exp(z_k)$, with $z_k = k\log(\delta/2) - b(p/2+k, \nu/2) - g(k+1)$, where $b(x,y) = \log B(x,y)$ is evaluated by `lbeta` in R.

To evaluate $p_*$, the main step is to obtain $u_* = \max\{u : u/\alpha = F(u)\}$, which can be rapidly approximated by the iteration $u_1 = 1$, $u_{n+1} = \alpha F(u_n)$. Then $p_* = F(u_*)$.

DEPARTMENT OF STATISTICS
UNIVERSITY OF CONNECTICUT
215 GLENBROOK ROAD, U-4120
STORRS, CONNECTICUT 06269
USA
E-MAIL: zchi@stat.uconn.edu